\documentclass[10pt,a4paper]{article}
\usepackage{graphics, graphicx}     

\usepackage{cite}
\usepackage{amsmath}
\usepackage{amssymb}
\usepackage{amsthm}
\theoremstyle{plain}
\newtheorem{Theorem}{Theorem}[section] %
\newtheorem{Lemma}{Lemma}[section]
\newtheorem{Proposition}{Proposition}[section]

\theoremstyle{definition}
\newtheorem{Remark}{Remark}[section]

\theoremstyle{definition}

\newtheorem{Problem}{Problem}[section]

\newenvironment{Proof} 
{\par\noindent{\it Proof of}} 
{\hfill$\vspace{5mm}\scriptstyle\blacksquare$} 

\numberwithin{equation}{section} 
\numberwithin{figure}{section} 
\numberwithin{table}{section} 

\begin{document}

\setcounter{page}{1}

\markboth{M.I. Isaev, R.G. Novikov}{Instability in the Gel'fand inverse problem
at high energies}

\title{Instability in the Gel'fand inverse problem
at high energies}
\date{}
\author{ { M.I. Isaev }}

\maketitle
\begin{abstract}
We give an instability estimate for the Gel'fand inverse boundary value problem at high energies.
Our instability estimate shows an optimality of several  important
preceeding stability results on inverse problems of such a type.
\end{abstract}

\section{Introduction}
In this paper we continue studies on the Gel'fand inverse boundary value  problem for  the  Schr\"odinger equation
\begin{equation}\label{eq} 
	-\Delta \psi  + v(x)\psi = E\psi, \ \  x \in  D,
\end{equation}
where
\begin{equation}\label{eq_c}
	\begin{aligned}
		 D &\text{ is an open bounded domain in } \mathbb{R}^d,\ d\geq 2, \ \\ &\text{with } \partial D \in C^2,
 	\end{aligned}
\end{equation}
\begin{equation}\label{eq_c2}
 v \in \mathbb{L}^{\infty}(D).
\end{equation}
As boundary data we consider the  map 
$\hat{\Phi} = \hat{\Phi}(E)$ such that
\begin{equation}\label{def_phi}
	\hat{\Phi}(E) (\psi|_{\partial D}) = \frac{\partial \psi}{\partial \nu}|_{\partial D}
\end{equation}
for all sufficiently regular solutions $\psi$ of (\ref{eq}) in $\bar{D} = D \cup \partial D,$ where $\nu$ is the outward normal to $\partial D$. Here we assume also that
\begin{equation}\label{correct}
	\text{$E$ is not a Dirichlet eigenvalue for  operator $-\Delta + v$ in $D$.} 
\end{equation}
The map $\hat{\Phi} = \hat{\Phi}(E)$ is known as the Dirichlet-to-Neumann map.

We consider the following inverse boundary value problem for equation (\ref{eq}): 

\begin{Problem}
 Given $\hat{\Phi}$ for some fixed $E$, find $v$. 
\end{Problem}

This problem is known as the Gel'fand inverse boundary value problem for the Schr\"odinger equation at fixed energy
(see \cite{Gelfand1954}, \cite{Novikov1988}). At zero energy this problem can be considered also as a generalization of the Calderon problem of the electrical impedance tomography (see \cite{Calderon1980}, \cite{Novikov1988}). 
Problem 1.1 can be also considered as an example of ill-posed problem: 
see \cite{LR1986}, \cite{BK2012} for an introduction to this theory.

There is a wide literature on the Gel'fand inverse problem at fixed energy. In a similar way with many other inverse problems, Problem 1.1 includes, in particular, the following questions: (a) uniqueness, (b) reconstruction,
(c) stability.

Global uniqueness results and global reconstruction methods for Problem 1.1 
were obtained for the first time in \cite{Novikov1988} in dimension $d \geq 3$  and in \cite{Buckhgeim2008} in dimension $d=2$.

Global logarithmic stability estimates for Problem 1.1 were obtained for the first time in \cite{Alessandrini1988} 
in dimension $d \geq 3$ and in \cite{NS2010} in dimension $d=2$.  
A principal improvement of the result of \cite{Alessandrini1988} was obtained recently in \cite{Novikov2011} (for the zero energy case):
stability of \cite{Novikov2011} optimally increases with increasing regularity of $v$. 

Note that for the Calderon problem (of the electrical impedance tomography) in its initial formulation the global uniqueness 
was firstly proved in \cite{SU1987} for $d\geq 3$ and in \cite{Nachman1996} for $d=2$.
Global logarithmic stability estimates for this problem were obtained
for the first time in \cite{Alessandrini1988} for $d\geq 3$
and \cite{Liu1997} for $d=2$. Principal increasing of global stability
of \cite{Alessandrini1988}, \cite{Liu1997} for the regular coefficient case was found in
\cite{Novikov2011} for $d\geq 3$ and \cite{S2011N} for $d=2$.
In addition, for the case of
piecewise real analytic conductivity the first uniqueness results for the Calderon problem in dimension $d\geq 2$
were given in 
\cite{KV1985}. Lipschitz stability estimate for the case of
piecewise constant conductivity was obtained in \cite{AV2005} (see \cite{Rondi2006} for additional studies in this direction).

The optimality of the logarithmic 
stability results of \cite{Alessandrini1988}, \cite{Liu1997} 
with their principal effectivizations of \cite{Novikov2011}, \cite{S2011N} (up to the value of the exponent)
follows from  \cite{Mandache2001}. 
An extention of the instability estimates of  \cite{Mandache2001}  to the case of the non-zero energy as well as to the case 
 of Dirichlet-to-Neumann map given on the energy intervals was obtained in \cite{IsaevDtN}.

On the other hand, it was found in \cite{Novikov1998}, \cite{Novikov2005+} (see also \cite{Novikov2008}, \cite{NS2012}) 
that for inverse problems for the  Schr\"odinger equation at fixed energy $E$ in dimension $d\geq 2$ (like Problem 1.1)
there is a H\"older stability modulo an error term rapidly decaying as $E \rightarrow +\infty$ (at least for the 
regular coefficient case). In addition, for Problem 1.1 for $d=3$, global energy dependent stability estimates changing
from logarithmic type to H\"older type for high energies were obtained in \cite{Isakov2011}, \cite{IN2012++}.
However, there is no efficient stability increasing with respect to increasing
coefficient regularity in the results of \cite{Isakov2011}.
 An additional study, motivated by \cite{Isakov2011}, \cite{Novikov2011}, was
given in \cite{NUW2011}.

The following stability estimate for Problem 1.1 was recently proved in \cite{IN2012++}:  
\begin{Theorem}[of \cite{IN2012++}]\label{Theorem_2.1} Let $D$ satisfy \eqref{eq_c}, where $d \geq 3$. 
Let  $v_j\in W^{m,1}(D)$, $m>d$, $supp\ v_j \subset D$ and $||v_{j}||_{W^{m,1}(D)} \leq N$  for some $N > 0$, $j = 1,2$, 
(where $W^{m,p}$ denotes the Sobolev space
of $m$-times smooth functions in $\mathbb{L}^p$). 
Let $v_1$, $v_2$ satisfy  \eqref{correct} 
for some fixed  $E\geq 0$.
Let $\hat{\Phi}_{1}(E)$ and $\hat{\Phi}_{2}(E)$ denote 
the DtN maps for
$v_1$ and  $v_2$, respectively. Let $s_1 = (m-d)/d$.   Then, for  any $\tau \in (0,1)$ and any $\alpha,\beta \in [0, s_1]$, $\alpha+\beta = s_1$, 
\begin{equation}\label{eq_2.1}
	||v_2 - v_1||_{L^\infty(D)} \leq A(1+\sqrt{E}) \delta^\tau + 
	B (1+\sqrt{E})^{-\alpha}\left(\ln\left(3+\delta ^{-1}\right)\right)^{-\beta}, 
\end{equation}
where $\delta  = ||\hat{\Phi}_{2}(E) - \hat{\Phi}_{1}(E) ||_{\mathbb{L}^\infty(\partial D) \rightarrow \mathbb{L}^\infty(\partial D)}$
and constants $A,B>0$ depend only on $N$, $D$, $m$, $\tau$.  
\end{Theorem}

In particular cases, H\"older-logarithmic stability estimate \eqref{eq_2.1} becomes coherent
(although less strong) with respect to results of  \cite{Novikov2005+}, \cite{Novikov2008}, \cite{Novikov2011}. 
In this connection we refer to \cite{IN2012++} for more detailed infromation.  Concerning two-dimensional analogs of results of Theorem \ref{Theorem_2.1}, see \cite{Novikov1998}, \cite{NS2012}, \cite{S2011N},
\cite{S2012}. 	

In a similar way with results of \cite{IN2012}, \cite{IN2012+}, estimate
	\eqref{eq_2.1} can be extended to the case when we do not assume that
	condition \eqref{correct} is fulfiled and consider an appropriate impedance boundary map (or Robin-to-Robin map)
	instead of the Dirichlet-to-Neumann map.

	In the present work we prove optimality of estimate \eqref{eq_2.1} (up to the values of the exponents $\alpha$, $\beta$) in dimension $d\geq 2$.  Our related instability results for 	Problem 1.1	are presented in Section 2, see Theorem \ref{Theorem_2.2} and Proposition \ref{Proposition_2.1}.  
	Their proofs are given in Section 4 and 
	are based on properties of solutions of the  Schr\"odinger equation in the unit ball given in Section 3.


\section{Main results}


In what follows we fix $D=B^d(0,1)$, where
\begin{equation}	
	B^d(x^0,\rho) = \{x \in \mathbb{R}^d: ||x-x^0||_{\mathbb{E}^d} < \rho\}, \ \ \ x_0 \in \mathbb{R}^d, \  \rho>0.
\end{equation}

Let 
\begin{equation}\label{norm_A}
	\begin{aligned}
		||F|| \text{ denote the norm of an operator}\\
		F: \mathbb{L}^{\infty}(\partial D) \rightarrow \mathbb{L}^{\infty}(\partial D).
	\end{aligned}
\end{equation}
We recall that if $v_1$, $v_2$ are potentials satisfying (\ref{eq_c2}), (\ref{correct})
for some fixed $E$, then
 \begin{equation}
	\hat{\Phi}_{2}(E) - \hat{\Phi}_{1}(E) \text { is a compact operator in } \mathbb{L}^{\infty}(\partial D),
\end{equation}
where $\hat{\Phi}_1$, $\hat{\Phi}_2$ are the DtN maps for $v_1$, $v_2$, respectively, see \cite{Novikov1988}, \cite{Novikov2005}.

Our main result is the following theorem:

\begin{Theorem}\label{Theorem_2.2}
Let $D = B^d(0,1)$, where $d\geq 2$.  Then for any fixed constants $A,B,\kappa, \tau, \varepsilon>0$, $m>d$ and $s_2>m$ there are 
some energy level $E>0$ and some potential $v \in C^m(D)$ such that condition 
\eqref{correct} holds for potentials $v$ and $v_0 \equiv 0$,  
simultaneously,  $\mbox{supp}\ v \subset D$, $\|v\|_{\mathbb{L}^\infty(D)} \leq \varepsilon$, $\|v\|_{C^m(D)} \leq  C_1$, where $C_1= C_1(d,m)>0$, but
\begin{equation}\label{eq_2.2}
	||v - v_0||_{L^\infty(D)} > A(1+\sqrt{E})^\kappa \delta^\tau + 
	B (1+\sqrt{E})^{2(s- s_2)}\left(\ln\left(3+\delta ^{-1}\right)\right)^{-s} 
\end{equation}
 for any $s\in[0,s_2]$, 
where
$\hat{\Phi}$, $\hat{\Phi}_{0}$ are the DtN map for $v$ and $v_0$, respectively, and
 $\delta  = ||\hat{\Phi}(E) - \hat{\Phi}_{0}(E) ||$ is defined according  to \eqref{norm_A}.
\end{Theorem}

 
 Theorem \ref{Theorem_2.2} shows, in particular, the optimality (at least for potentials in the neighborhood of zero) of 
 estimate \eqref{eq_2.1} 
(up to the values of the exponents $\alpha$, $\beta$). 
As a corollary of Theorem \ref{Theorem_2.2}, one can obtain an optimality 
of the stability results of \cite{Novikov1998}, \cite{Novikov2005+}, \cite{Novikov2008}, \cite{NS2012}.

In the present work Theorem \ref{Theorem_2.2} is proved by explicit instability example with complex
potentials. Examples of this type were considered for the 
first time in \cite{Mandache2001} for showing the exponential instability in Problem 1.1 
in the zero energy case. An extention to the case of the non-zero energy as well as to the case  of Dirichlet-to-Neumann map given on the energy intervals was obtained in \cite{IsaevDtN}.

Let us consider the cylindrical variables:
\begin{equation}
	\begin{aligned}
		(r_1,\theta,x')\in \mathbb{R}_+ \times \mathbb{R}/2\pi\mathbb{Z} \times \mathbb{R}^{d-2}, \\
		 r_1\cos\theta = x_1, \ \ r_1 \sin\theta =x_2, \\x' = (x_3,\ldots,x_d).
	\end{aligned}
\end{equation}
Take $\phi\in C^\infty(\mathbb{R}^2)$ with support in $B^2(0,1/3) \cap \{x_1>1/4\}$ and with 
$\|\phi\|_{\mathbb{L}^\infty} =1$.
For integers $m,n>0$, define the complex potential
	\begin{equation}\label{v_nm}
		v_{nm} = n^{-m}e^{in\theta} \phi(r_1, |x'|).
	\end{equation}
We recall that
\begin{equation}\label{eq_v_nm}
	\|v_{nm}\|_{\mathbb{L}^\infty} = n^{-m}, \ \ \ \|v_{nm}\|_{C^m} \leq C_1, 
\end{equation}
where $C_1=C_1(d,m)>0$. 
Note that $C_1$ is the same as in Theorem \ref{Theorem_2.2}.
Estimates \eqref{eq_v_nm} were given in \cite{Mandache2001} (see Theorem 2 of \cite{Mandache2001}).

To prove Theorem \ref{Theorem_2.2} we use, in partucular,    the following proposition:
\begin{Proposition}\label{Proposition_2.1}
Let $D = B^d(0,1)$, where $d\geq 2$.  Let condition \eqref{correct} hold with 
  $v\equiv v_{nm}$  (of \eqref{v_nm}) and $v\equiv v_0\equiv 0$ 
for some $E>0$ and some integers $m>0$, $n>20(1+\sqrt{E})^2$. Then, for any $\sigma>0$,
\begin{equation}\label{eq_Pr_2.1}
	\begin{aligned}
	\|\hat{\Phi}_{nm}(E) - \hat{\Phi}_0(E)\|_{H^{-\sigma}(\mathbb{S}^{d-1})\rightarrow H^{\sigma}(\mathbb{S}^{d-1})}
		\leq C_2(1+Q+EQ) 2^{-n/4},
	\end{aligned}
\end{equation}
where  $\hat{\Phi}_{nm}$, $\hat{\Phi}_{0}$ are the DtN map for $v_{nm}$ and $v_0$, respectively, 
$C_2 = C_2(d,\sigma)>0$, 
\begin{equation}
		Q = \|(-\Delta +v_0 - E)^{-1}\|_{\mathbb{L}^2(D) \rightarrow \mathbb{L}^2(D)} +
		\|(-\Delta + v_{nm} -E)^{-1}\|_{\mathbb{L}^2(D) \rightarrow \mathbb{L}^2(D)},
\end{equation}
where $(-\Delta+v_0- E)^{-1}$, 
$(-\Delta + v_{nm} -E)^{-1}$ are
considered with the Dirichlet boundary condition in $D$ and $H^{\pm\sigma} = W^{\pm\sigma,2}$ denote the standart Sobolev spaces.
\end{Proposition}

Analogs of  estimate \eqref{eq_Pr_2.1} (but without dependence of the energy) 
were given in Theorem 2 of \cite{Mandache2001} for the zero energy case  and in
Theorem 2.4 of \cite{IsaevDtN} for the case of the non-zero energy and the case  of the energy intervals.

We obtain Theorem \ref{Theorem_2.2}, combining known results on the spectrum of the Laplace operator in the unit ball 
(see formula \eqref{Weyl} below), 
Proposition 2.1, estimates \eqref{eq_v_nm} and the fact  that 
 	\begin{equation}\label{HsigmaL}
 		\|F\|_{L^{\infty}(\mathbb{S}^{d-1}) \rightarrow L^{\infty}(\mathbb{S}^{d-1})} 
 		\leq c(d,\sigma)\|F\|_{H^{-\sigma}(\mathbb{S}^{d-1}) \rightarrow H^{\sigma}(\mathbb{S}^{d-1})}
 	\end{equation}
for sufficiently large $\sigma$. The detailed proof of  Theorem \ref{Theorem_2.2} 
and the proof of Proposition \ref{Proposition_2.1}
are given in Section 4. These proofs use, in particular, results, presented in Section 3.

\begin{Remark}
 	In a similar way with \cite{Mandache2001}, \cite{IsaevDtN}, using  a ball packing and covering by ball arguments 
 	(see also \cite{CR2003}),
 	the instability result of Theorem \ref{Theorem_2.2} can be extended to the case when 
 	only real-valued potentials  are considered and in the neighborhood of any potential (not only $v_0\equiv 0$).
\end{Remark}

\section{Some  properties of solutions of the  Schr\"odinger equation in the unit ball}
In this section we continue assume that $D = B^d(0,1)$, where $d\geq 2$. 
We fix an orthonormal basis in $\mathbb{L}^2(\mathbb{S}^{d-1}) = \mathbb{L}^2(\partial D)$
\begin{equation}
	\begin{array}{l}
	\displaystyle
		\{f_{jp} : j \geq 0,\  1 \leq p \leq p_j \}, \\
		\text{$f_{jp}$ is a spherical harmonic of degree $j$,}
	\end{array}
\end{equation}
   where $p_j$  is the dimension of the space of spherical harmonics of order $j$,
\begin{equation}
	p_j = \binom {j+d-1} {d-1} - \binom {j+d-3} {d-1},
\end{equation}
	where 
\begin{equation}
	\binom {n} {k} = 
		\frac{n(n-1)\cdots(n-k+1)}{k!} \ \ \ \text{ for $n \geq 0$} 
\end{equation}	
and
\begin{equation}
	\binom {n} {k} = 
		0 \ \ \ \text{ for $n < 0$.} 
\end{equation}
The precise choice of $f_{jp}$ is irrelevant for our purposes. Besides orthonormality, we only need
$f_{jp}$ to be the restriction of a homogeneous harmonic polynomial of degree $j$ to the sphere
and so $|x|^j f_{jp}(x/|x|)$ is harmonic. 
We use also the polar coordinates 
$(r, \omega) \in \mathbb{R}_+ \times \mathbb{S}^{d-1}$, with $x = r\omega \in \mathbb{R}^d$.

\begin{Lemma}\label{Lemma_3.1}
	Let $D = B^d(0,1)$, where $d \geq 2$. 
Let  potential $v$ satisfy  \eqref{eq_c2} and \eqref{correct} 
for some fixed  $E$. Let $||v||_{\mathbb{L}^\infty(D)} \leq N$, for some $N > 0$. Then 
for any solution $\psi\in C(D\cup\partial D)$ of equation \eqref{eq} the following inequality holds:
\begin{equation}\label{3.6}
	\begin{aligned}
	\|\psi\|_{\mathbb{L}^2(D)} \leq \Big(1+(N+|E|) 
	\|(-\Delta + v -E)^{-1}\|_{\mathbb{L}^2(D) \rightarrow \mathbb{L}^2(D)}\Big) \|f\|_{\mathbb{L}^2(\partial D)}, 
	\end{aligned}
\end{equation}
where $f = \psi|_{\partial D}$,  $(-\Delta + v -E)^{-1}$ is
considered with the Dirichlet boundary condition in $D$.
\end{Lemma}
\begin{Proof} {\it Lemma \ref{Lemma_3.1}.}
We expand the function $f$ in the basis $\{f_{jp}\}$: 
\begin{equation}
	f = \sum\limits_{j,p} c_{jp} f_{jp}. 
\end{equation}
We have that
\begin{equation}\label{3.8}
	 \|f\|^2_{\mathbb{L}^2(\partial D)} = \sum\limits_{j,p} |c_{jp}|^2.
\end{equation}
Let 
\begin{equation}
	\psi_0(x) = \sum\limits_{j,p} c_{jp} r^{j}f_{jp}(\omega).
\end{equation}
Note that
\begin{equation}\label{3.10}
	\begin{aligned}
	\|\psi_0\|^2_{\mathbb{L}^2(D)}  = \sum\limits_{j,p} |c_{jp}|^2 \|r^{j}f_{jp}(\omega)\|^2_{\mathbb{L}^2(D)}	
		=\\= \sum\limits_{j,p} |c_{jp}|^2 \int_0^1 r^{2j+d-1} dr \leq  \sum\limits_{j,p} |c_{jp}|^2
	\end{aligned}
\end{equation}
 Using \eqref{eq} and the fact
that $\psi_0$ is harmonic, we get that
\begin{equation}\label{3.11}
	(-\Delta+v - E) (\psi - \psi_0) = (E-v)\psi_0.
\end{equation}
Since $\psi|_{\partial D} = \psi_0|_{\partial D} = f$, using \eqref{3.11}, we find that
\begin{equation}\label{3.12}
	\|\psi - \psi_0 \|_{\mathbb{L}^2(\partial D)} \leq 
	(N+|E|) 	\|(-\Delta + v -E)^{-1}\|_{\mathbb{L}^2(D) \rightarrow \mathbb{L}^2(D)} 
	\|\psi_0\|_{\mathbb{L}^2(D)}.
\end{equation}
Combining \eqref{3.8}, \eqref{3.10}, \eqref{3.12}, we obtain \eqref{3.6}.
\end{Proof}
  
Let $<\cdot,\cdot>$ denote the scalar product in the Hilbert space $\mathbb{L}^2(\partial D)$:
\begin{equation}
	<f,g> = \int\limits_{\partial D} f(x) \bar{g}(x) dx,  
\end{equation}
where $f,g \in \mathbb{L}^2(\partial D)$.
  
\begin{Lemma}\label{Lemma_3.2}
	Let $D = B^d(0,1)$, where $d \geq 2$. 
	Let  potentials $v_1$, $v_2$ satisfy  \eqref{eq_c2} and \eqref{correct} 
		for some fixed  $E$. Let $v_1$, $v_2$ be supported in $B^d(0,1/3)$ and $||v_i||_{\mathbb{L}^\infty(D)} \leq N$, 
		$i=1,2$,  for some $N > 0$. Then for any 
		$j_1,j_2 \in \mathbb{N}\cup\{0\}$,  $1 \leq p_1 \leq p_{j_1}$, $1 \leq p_2 \leq p_{j_2}$ and 
		$j_{max}= \max\{j_1,j_2\} \geq 10(1+\sqrt{|E|})^2$ the following inequality holds:
		\begin{equation}\label{eq_3.2}
			\left|\left\langle f_{j_1 p_1}, \left(\hat{\Phi}_1(E) - \hat{\Phi}_2(E)\right)f_{j_2 p_2}\right\rangle\right| 
			\leq C(d)\Big(1+(N+|E|)Q\Big)2^{-j_{max}},
		\end{equation}
		where
		\begin{equation}
				Q = \|(-\Delta +v_1 - E)^{-1}\|_{\mathbb{L}^2(D) \rightarrow \mathbb{L}^2(D)} +
				\|(-\Delta + v_2 -E)^{-1}\|_{\mathbb{L}^2(D) \rightarrow \mathbb{L}^2(D)},
		\end{equation} 
$\hat{\Phi}_1$, $\hat{\Phi}_2$ are the DtN map for $v_1$ and $v_2$, respectively, and $(-\Delta+v_1- E)^{-1}$, 
$(-\Delta + v_2 -E)^{-1}$ are
considered with the Dirichlet boundary condition in $D$.
\end{Lemma}

Analogs of  estimate \eqref{eq_3.2} (but without dependence of the energy) 
were given in Lemma 1 of \cite{Mandache2001} for the zero energy case  and in
Lemma 3.4 of \cite{IsaevDtN} for the case of the non-zero energy and the case  of the energy intervals. 


  We prove Lemma \ref{Lemma_3.2} for $E\neq 0$ in Section 5, 
  using expression of solutions of equation $-\Delta\psi = E\psi$ 
in $B^d(0,1)\setminus B^d(0,1/3)$ in terms of the Bessel functions $J_\alpha$ and $Y_\alpha$ with integer or half-integer
order $\alpha$.

\section{Proofs of Proposition 2.1 and Theorem \ref{Theorem_2.2}}
We continue to assume that $D = B^d(0,1)$, where $d\geq 2$ and to use the 
orthonormal basis 
		$\{f_{jp} : j\in\mathbb{N}\cup\{0\} ,\  1 \leq p \leq p_j \}$ 
		in $\mathbb{L}^2(\mathbb{S}^{d-1}) = \mathbb{L}^2(\partial D)$. 
The Sobolev spaces $H^\sigma (\mathbb{S}^{d-1})$ can be defined by
 	\begin{equation}
 		\begin{aligned}
 			\left\{ \sum_{j,p} c_{jp}f_{jp}: \Big\| \sum_{j,p} c_{jp}f_{jp}\Big\|_{H^\sigma}<+\infty\right\},\\
			\Big\| \sum_{j,p} c_{jp}f_{jp}\Big\|^2_{H^\sigma} = \sum_{j,p}(1+j)^{2\sigma}|c_{jp}|^2,
		\end{aligned}
	\end{equation} 
see, for example, \cite{Mandache2001}.
	
Consider an operator $A: H^{-\sigma}(\mathbb{S}^{d-1}) \rightarrow H^\sigma(\mathbb{S}^{d-1})$.	
We denote its matrix elements in the basis $\{f_{jp}\}$ by 
\begin{equation}
a_{j_1 p_1 j_2 p_2} = <f_{j_1 p_1},A f_{j_2 p_2}>.
\end{equation}
We identify in the sequel an operator $A$ with its matrix 
$\{a_{j_1 p_1 j_2 p_2}\}$. In this section we always assume that  
$j_1, j_2 \in \mathbb{N}\cup\{0\}$,   $1 \leq p_1 \leq p_{j_1}$,  $1 \leq p_2 \leq p_{j_2}$.

We recall that (see formula (12) of \cite{Mandache2001})
\begin{equation}\label{AHS}
	\|A\|_{H^{-\sigma}(\mathbb{S}^{d-1}) \rightarrow H^\sigma(\mathbb{S}^{d-1})} 
	\leq 4 \sup\limits_{j_1,p_1,j_2,p_2} (1+\max\{j_1,j_2\})^{2\sigma+d}|a_{j_1 p_1 j_2 p_2}|.
\end{equation}

\begin{Proof} {\it Proposition 2.1.}
In a similar way with the proof of Theorem 2 of \cite{Mandache2001} we obtain that
	\begin{equation}\label{jizero}
		<f_{j_1p_1}, \left(\hat{\Phi}_{mn}(E) - \hat{\Phi}_0(E)\right)f_{j_2p_2}> = 0
	\end{equation}  
	for $j_{max} = \max\{j_1,j_2\} \leq \left[\frac{n-1}{2}\right]$ (the only difference is that instead of the operator $-\Delta$ we consider 	the operator $-\Delta - E$), where $[ \cdot ]$ denotes the integer part of a number. 
	Note that 
	\begin{equation}\label{note_2.1}
		\left[\frac{n-1}{2}\right]+1 \geq n/2 > 10(1+\sqrt{E})^2, \ \ \ \|v_{nm}\|_{\mathbb{L}^\infty(D)}\leq 1.
	\end{equation}
		Combining (\ref{AHS}),  (\ref{jizero}), \eqref{note_2.1} and Lemma \ref{Lemma_3.2}, we get that
	\begin{equation}
		\begin{aligned}
	 \|\hat{\Phi}_{mn}(E) - &\hat{\Phi}_0(E)\|_{H^{-\sigma}(\mathbb{S}^{d-1}) \rightarrow H^\sigma(\mathbb{S}^{d-1})} \leq\\
	 &\leq 4C(d)\Big(1+(1+E)Q\Big)\sup_{j_{max}\geq n/2}(1+j_{max})^{2\sigma+d}2^{-j_{max}} \leq \\
	 &\leq C_2(d,\sigma) (1+Q+EQ) 2^{-n/4},
	 \end{aligned}
	\end{equation}
	where 
\begin{equation}
		Q = \|(-\Delta +v_0 - E)^{-1}\|_{\mathbb{L}^2(D) \rightarrow \mathbb{L}^2(D)} +
		\|(-\Delta + v_{nm} -E)^{-1}\|_{\mathbb{L}^2(D) \rightarrow \mathbb{L}^2(D)}.
\end{equation}
\end{Proof}

Let $N(\rho)$ denote the counting function of the Laplace operator in $D$
\begin{equation}
	N(\rho) = |\{\lambda<\rho^2: \lambda \text{ is a Dirichlet eigenvalue of  $-\Delta$  in } D \}|, 
\end{equation}
where $|\cdot|$ is the cardinality of the corresponding set.
We recall that according to the Weyl formula (of \cite{Weyl1912}):
\begin{equation}\label{Weyl}
N(\rho)\leq c_1(d)\rho^d.
\end{equation}

\begin{Lemma} Let $D=B^d(0,1)$, where $d\geq 1$. Then
	for any $\rho>1$ there is some $E=E(\rho)\in (\rho^2, 2\rho^2)$ such that the interval
	\begin{equation}\label{Pr_4.1}
		\left(E(\rho)-c_2\rho^{2-d}, E(\rho)+c_2\rho^{2-d}\right) 
	\end{equation}
		 does not contain Dirichlet eigenvalues of $-\Delta$  in  $D$, where $c_2=c_2(d)>0$.
\end{Lemma}
\begin{Proof} {\it Lemma 4.1.}
We put $c_2 = 2^{d-1}/(c_1(d)+1)$. Then  we can select 
$k$ disjoint intervals of the length $2c_2\rho^{2-d}$  in the interval $(\rho^2, 2\rho^2)$, where 
	\begin{equation}
		k = \left[\frac{\rho^2}{2c_2\rho^{2-d}}\right] = [(c_1(d)+1)\rho^d] > N(\rho).
	\end{equation}
Thus, we have that at least one of these intervals  does not contain Dirichlet eigenvalues of $-\Delta$  in  
$D = B^d(0,1)$.
\end{Proof}

\begin{Proof} {\it Theorem \ref{Theorem_2.2}.}
 Let $E = E(\rho)$ be the number of Lemma 4.1 for some $\rho>1$. Using \eqref{Pr_4.1}, we find that
  the distance from $E$ to the  Dirichlet spectrum of
  the operator $-\Delta$  in  $D$ is not less than   $c_2\rho^{2-d}$. Using also that $E\in(\rho^2,2\rho^2)$, we get that
 \begin{equation}\label{4.12}
 	\|(-\Delta  - E)^{-1}\|_{\mathbb{L}^2(D) \rightarrow \mathbb{L}^2(D)} \leq \frac{1}{c_2\rho^{2-d}} \leq E^{(d-2)/2}/c_2, 
 \end{equation}
where $(-\Delta  -E)^{-1}$ is
considered with the Dirichlet boundary condition in $D$.
 Let  
 \begin{equation}\label{4.13}
 	n = [20(1+\sqrt{E})^2] + 1. 
 \end{equation}
 Using  \eqref{eq_v_nm} and \eqref{Pr_4.1}, we find that
  the distance from $E$ to the  Dirichlet spectrum of
  the operator $-\Delta+v_{nm}$  in  $D$ is not less than   $c_2\rho^{2-d} - n^{-m}$,  where 
  $v_{nm}$ is defined according to \eqref{v_nm}. Since $m>d$ and $E\in(\rho^2,2\rho^2)$, 
  using \eqref{4.13}, we get that
 \begin{equation}\label{4.14}
 	\begin{aligned}
 	 	\|(-\Delta + v_{nm} - E)^{-1}\|_{\mathbb{L}^2(D) \rightarrow \mathbb{L}^2(D)} \leq c_3 E^{(d-2)/2},\\
 	 		E=E(\rho), \ \ \ \rho \geq \rho_1(d,m)>1,\\
 	 		c_3 = c_3(d,m)>0,
 	 \end{aligned}	
 \end{equation}
  where $(-\Delta + v_{nm} -E)^{-1}$ is
considered with the Dirichlet boundary condition in $D$.
  
 Combining Proposition \ref{Proposition_2.1} and estimates \eqref{HsigmaL}, \eqref{4.12}, \eqref{4.14}, we find that
 \begin{equation}\label{4.15}
 	\begin{aligned}
 		\delta = \|\hat{\Phi}_{nm}(E) - \hat{\Phi}_0(E)\|_{
 		\mathbb{L}^{\infty}(\mathbb{S}^{d-1})\rightarrow \mathbb{L}^{\infty}(\mathbb{S}^{d-1})}
		\leq c_4 E^{d/2} 2^{-n/4},\\
		 E=E(\rho), \ \ \ \rho \geq \rho_1(d,m)>1,\\
		  	n = [20(1+\sqrt{E})^2] + 1\\
		  	c_4 = c_4(d,m)>0.
	\end{aligned}
 \end{equation} 
 Since $s_2>m$, taking $\rho$ big enough and using \eqref{4.15}, we obtain the following inequalities:
 \begin{equation}\label{4.16}
 	 n^{-m} < \varepsilon,
 \end{equation}
 \begin{equation}
 	A(1+\sqrt{E})^\kappa \delta^\tau < \frac{1}{2}n^{-m},
 \end{equation}
 \begin{equation}
 	\begin{aligned}
 	B (1+\sqrt{E})^{2(s- s_2)}\left(\ln\left(3+\delta ^{-1}\right)\right)^{-s} &<
 	\frac{1}{2}n^{-m}, 	
 	\\
 		0 \leq s &\leq s_2,
 	\end{aligned}
 \end{equation}
 where 
 \begin{equation}\label{4.19}
 	E=E(\rho), \ \ \ 	n = [20(1+\sqrt{E})^2] + 1.
 \end{equation}
 Combining   \eqref{v_nm}, \eqref{eq_v_nm}, \eqref{4.16} - \eqref{4.19}, we get that
 \begin{equation}
 	\begin{aligned}
 	A(1+\sqrt{E})^\kappa \delta^\tau + B (1+\sqrt{E})^{2(s- s_2)}\left(\ln\left(3+\delta ^{-1}\right)\right)^{-s} <
 	\\ <\frac{1}{2}n^{-m} + \frac{1}{2}n^{-m} = 
 \|v_{nm}- v_0\|_{\mathbb{L}^{\infty}(D)} \\
		\|v_{nm}\|_{\mathbb{L}^{\infty}(D)} = n^{-m} < \varepsilon, \\  	\|v_{nm}\|_{C^{m}(D)} < C_1, \\   \mbox{supp}\, v_{nm} \subset D.
	\end{aligned}
 \end{equation}
\end{Proof}

\section{Proof of Lemma \ref{Lemma_3.2}}
To prove Lemma \ref{Lemma_3.2} we need some preliminaries.
Consider the problem of finding solutions of the form $\psi(r,\omega) = R(r)f_{jp}(\omega)$
of equation (\ref{eq}) with $v\equiv 0$ and $D = B^d(0,1)$, where $d\geq 2$. We recall that:
\begin{equation}
	\Delta = \frac{\partial^2}{(\partial r)^2} + (d-1)r^{-1} \frac{\partial}{\partial r} + r^{-2}\Delta_{S^{d-1}},   
\end{equation}
where $\Delta_{S^{d-1}}$ is  Laplace-Beltrami operator on $S^{d-1}$,
\begin{equation}
	\Delta_{S^{d-1}}f_{jp} = -j(j+d-2)f_{jp}.
\end{equation}
Then we obtain the following equation for $R(r)$:
\begin{equation}
	-R'' - \frac{d-1}{r}R' + \frac{j(j+d-2)}{r^2}R = ER.
\end{equation}
Taking $R(r) = r^{-\frac{d-2}{2}}\tilde{R}(r)$, we get
\begin{equation}
	r^2\tilde{R}'' + r\tilde{R}' + \left(Er^2 - \left(j+\frac{d-2}{2}\right)^2\right)\tilde{R} = 0.
\end{equation}
This equation is known as the Bessel equation. For $E =k^2\neq 0$ it has two linearly independent solutions 
$J_{j+\frac{d-2}{2}}(kr)$ and $Y_{j+\frac{d-2}{2}}(kr),$ where 
\begin{equation}\label{jalphaz}
	J_\alpha(z) = \sum\limits_{m=0}\limits^{\infty}
		\frac{(-1)^m(z/2)^{2m+\alpha}}{\Gamma(m+1)\Gamma(m+\alpha + 1)},
\end{equation} 
\begin{equation}
	Y_\alpha(z) = \frac{J_\alpha(z)\cos\pi\alpha  - J_{-\alpha}(z)}{\sin\pi\alpha} \text{ for  $\alpha \notin \mathbb{Z}$,}
\end{equation}
and 
\begin{equation}
	Y_\alpha(z) = \lim\limits_{\alpha'\rightarrow \alpha}Y_{\alpha'}(z) \text{ for  $\alpha \in \mathbb{Z}$.}
\end{equation}
We recall also that the system of functions 
		\begin{equation}\label{sys1}
			\begin{aligned}
				\left\{\psi_{jp}(r,\omega) = R_j(k,r)f_{jp}(\omega) : j \in \mathbb{N}\cup\{0\}, 1\leq p\leq p_j \right\}, \\
	\text{
	is complete orthogonal system (in the sense of $\mathbb{L}^2$) in the space}\\ \text{of solutions of equation (\ref{eq}) 
	in 
	$D' = B(0,1)\setminus B(0,1/3)$ } \\
	\text{with $v \equiv 0$, $E = k^2$ and boundary condition $\psi|_{r=1} = 0$},
	\end{aligned}
\end{equation}
		where 
		\begin{equation}\label{sys2}
			R_j(k,r) = {r^{-\frac{d-2}{2}}} \Big(Y_{j+\frac{d-2}{2}}(kr)J_{j+\frac{d-2}{2}}(k) - 
			J_{j+\frac{d-2}{2}}(kr)Y_{j+\frac{d-2}{2}}(k) \Big).
		\end{equation}	
	 For the proof of \eqref{sys1} see, for example, \cite{IsaevDtN}.

\begin{Lemma}\label{Bessel}
	For any $\rho>0$, integers $d \geq 2$,  $n\geq 10(\rho+1)^2$ and  $z\in\mathbb{C}$, $|z|\leq \rho$, the
	following inequalities hold:
		\begin{equation}\label{Jn}
			\frac{1}{2}\frac{(|z|/2)^\alpha}{\Gamma(\alpha+1)} \leq |J_\alpha(z)| \leq \frac{3}{2}\frac{(|z|/2)^\alpha}{\Gamma(\alpha+1)},
		\end{equation}
		\begin{equation}\label{Jn'}
			|J'_\alpha(z)| \leq {3}\frac{(|z|/2)^{\alpha-1}}{\Gamma(\alpha)},
		\end{equation}
		\begin{equation}\label{Yn}
			\frac{1}{2\pi}(|z|/2)^{-\alpha}\Gamma(\alpha) \leq |Y_\alpha(z)| \leq \frac{3}{2\pi}(|z|/2)^{-\alpha}\Gamma(\alpha)
		\end{equation}
		\begin{equation}\label{Yn'}
			|Y'_\alpha(z)| \leq \frac{3}{\pi}(|z|/2)^{-\alpha-1}\Gamma(\alpha+1)
		\end{equation}
		where $'$ denotes derivation with respect to $z$, $\alpha = n+\frac{d-2}{2}$ and $\Gamma(x)$ is the Gamma function.
\end{Lemma}  
In fact, the proof of Lemma \ref{Bessel} is given in \cite{IsaevDtN} (see Lemma 3.3  of \cite{IsaevDtN}).
It was shown in \cite{IsaevDtN} that inequalities \eqref{Jn} - \eqref{Yn'} hold for any $n>n_0$, where $n_0$ is such that
\begin{equation}\label{5.14}
\left\{
		\begin{array}{l}\displaystyle
				n_0>3,
			\\\displaystyle
			\exp\left(\frac{\rho^2/{4}}{n_0+1}\right) - 1 \leq 1/2 ,
				\\\displaystyle
				3\pi\frac{\max\left(1,(\rho/2)^{2n_0+1}\right)}{\Gamma(n_0)} + 
			\frac{\rho^2}{2n_0-\rho^2} + \frac{(\rho/2)^{2n_0}e^{\rho^2/4}}{\Gamma(n_0)} \leq 1/2,
		\end{array}\right.
	\end{equation}
(see  formula (6.18) of \cite{IsaevDtN}).
The only thing to check is that $n_0 = [10(\rho+1)^2]-1$ satisfy \eqref{5.14}, where $[ \cdot ]$ denotes the integer part of a number, 
The first two inequalities are obvious. The third follows from the estimate 
\begin{equation}
\Gamma(n_0) = (n_0-1)!\geq \left(\frac{n_0-1}{e}\right)^{n_0-1}.
\end{equation}	

The final part of the proof of Lemma \ref{Lemma_3.2} consists of the following:
	first, we consider the case when $E = k^2 \neq 0$ and 
	\begin{equation}\label{as_j}	
		j_1 = \max\{j_1,j_2\}\geq 10(1+|k|)^2.
	\end{equation}		 
	Let $\psi_1$, $\psi_2$ denote the solutions of equation (\ref{eq})   with 
	boundary condition $\psi|_{\partial D} = f_{j_2p_2}$ and potentials $v_1$ and $v_2$, respectively. 
	Using Lemma 3.1 for  $v_1$ and $v_2$, we get that
	\begin{equation}\label{l2sol}
		\|\psi_1 - \psi_2\|_{\mathbb{L}^2(D)}\leq 2\Big(1+(N+|E|)Q\Big),
	\end{equation}
	where
		\begin{equation}
				Q = \|(-\Delta +v_1 - E)^{-1}\|_{\mathbb{L}^2(D) \rightarrow \mathbb{L}^2(D)} +
				\|(-\Delta + v_2 -E)^{-1}\|_{\mathbb{L}^2(D) \rightarrow \mathbb{L}^2(D)},
		\end{equation} 
	
		Note that $\psi_1 - \psi_2$ is the solution of equation (\ref{eq}) in 
	$D' = B(0,1)\setminus B(0,1/3)$ with potential $v \equiv 0$ and boundary condition $\psi|_{r=1} = 0.$
		According to \eqref{sys1}, we have that
	\begin{equation}\label{ck}
		\psi_1 - \psi_2 = \sum\limits_{j,p} c_{jp}\psi_{jp} \  \text{ in $D'$} 
	\end{equation}
	for some $c_{jp}$, where
	\begin{equation}
		\psi_{jp}(r,\omega) = R_{j}(k,r)f_{jp}(\omega).
	\end{equation}  
	Since $R_j(k,1) = 0$, we find that
	\begin{equation}
		\left.\frac{\partial R_j(k,r)}{\partial r}\right|_{r=1} = 
				\left.\frac{\partial \left({r^{\frac{d-2}{2}}}R_j(k,r)\right)}{\partial r}\right|_{r=1}.
	\end{equation}
	For $j \geq 10(1+|k|)^2$, using Lemma \ref{Bessel}, we have that	
		\begin{equation}\label{R'+}
			\begin{array}{c}
			\displaystyle
				\left|
				\frac{ \left.\frac{\partial R_i(k,r)}{\partial r}\right|_{r=1} }{Y_{\alpha}(k)J_{\alpha}(k)}  
				\right|
				=
				|k|
				\left|
					\frac{Y'_{\alpha}(k)}{Y_{\alpha}(k)}
					-
					\frac{J'_{\alpha}(k)}{J_{\alpha}(k)}
				\right| \leq
			\\\displaystyle
				\leq
				6 |k|
				\left(
						\frac
						{(|k|/2)^{-\alpha-1}\Gamma(\alpha+1)}   
						{(|k|/2)^{-\alpha}\Gamma(\alpha)}  
				+
						\frac
						{(|k|/2)^{\alpha-1}\Gamma(\alpha+1)}   
						{(|k|/2)^{\alpha}\Gamma(\alpha)} 
				\right) 
				=
				6\alpha,
			\end{array}
		\end{equation}
		\begin{equation}\label{ylemma1}
			\begin{aligned}
			\bigg(
			\frac
			{||r^{-\frac{d-2}{2}}Y_\alpha(kr)||_{\mathbb{L}^2(\{1/3 <|x|<2/5\})}}
			{|Y_{\alpha}(k)|} 
			\bigg)^2 
			&\geq \\ \geq
			\int_{1/3}^{2/5} 
				\left(
						\frac{1}{3}
						\frac
						{(|k|r/2)^{-\alpha}\Gamma(\alpha)}   
						{(|k|/2)^{-\alpha}\Gamma(\alpha)}  
				\right)^2 r\, dr 
				&\geq 
				\left(\frac{2}{5} - \frac{1}{3}\right)\frac{1}{3} \left(\frac{1}{3} (5/2)^\alpha \right)^{2},
			\end{aligned}
		\end{equation}
		\begin{equation}\label{jlemma1}
			\begin{aligned}
			\left(
			\frac
			{||r^{-\frac{d-2}{2}}J_\alpha(kr)||_{\mathbb{L}^2(\{1/3 <|x|<2/5\})}}
			{|J_{\alpha}(k)|} 
			\right)^2
			&\leq \\ \leq
			\int_{1/3}^{2/5} 
				\left(
						3
						\frac
						{(|k|r/2)^{\alpha}\Gamma(\alpha)}   
						{(|k|/2)^{\alpha}\Gamma(\alpha)}  
				\right)^2 r\, dr 
				&\leq 
				\left(\frac{2}{5} - \frac{1}{3}\right)\frac{1}{3} \left(3 (2/5)^\alpha\right)^{2},
			\end{aligned}
		\end{equation}
		where $\alpha = j + \frac{d-2}{2}.$		Note that if $j \geq 10(1+|k|)^2$ then $j + \frac{d-2}{2} > 3$. 
		Combining (\ref{ylemma1}) and (\ref{jlemma1}), we get that
	\begin{equation}\label{psi+}
		\begin{aligned}
		\frac
		{||\psi_{jp}||_{L^2(\{1/3 <|x|<2/5\})}}
		{\left|Y_{\alpha}(k)J_{\alpha}(k)\right|}
		&\geq \\ \geq
		\left(\Big(\frac{2}{5} - \frac{1}{3}\Big)\frac{1}{3}\right)^{1/2}
		&\left(\frac{1}{3}(5/2)^\alpha - 3 (2/5)^\alpha\right) 
		> \frac{6}{1000}(5/2)^{\alpha}.
		\end{aligned}
\end{equation}
	Combining \eqref{R'+} and \eqref{psi+}, we find that
	\begin{equation}\label{Lambda1}
		\left|\left.\frac{\partial R_j(k,r)}{\partial r}\right|_{r=1}\right| \leq 1000\alpha(5/2)^{-\alpha}
		||\psi_{jp}(E)||_{\mathbb{L}^2(\{1/3 <|x|<1\})}.
	\end{equation}		
		Proceeding from (\ref{ck}) and using the Cauchy-Schwarz inequality, we get that
	\begin{equation}\label{Lambda2}
		|c_{jp}| = 
		\left|
		\frac
		{\Big< \psi_{jp}, \psi_1 - \psi_2 \Big>_{\mathbb{L}^2(\{1/3 <|x|<1\})}}
		{||\psi_{jp}(E)||^2_{\mathbb{L}^2(\{1/3 <|x|<1\})}} 
		\right|
		\leq
		\frac
		{||\psi(E) - \psi_0(E)||_{\mathbb{L}^2(B(0,1))}}
		{||\psi_{jp}(E)||_{\mathbb{L}^2(\{1/3 <|x|<1\})}}.
	\end{equation}
	Using \eqref{ck}, we find that
	\begin{equation}\label{5.28}
		\begin{aligned}
		\left\langle  f_{j_1p_1}, \left(\hat{\Phi}_1(E) - \hat{\Phi}_2(E)\right) f_{j_2p_2}\right\rangle
		=
		\left\langle 
		f_{j_1p_1},  
		\left.\frac
		{\partial(\psi_1-\psi_2)}
		{\partial\nu}\right|_{\partial D}
		\right\rangle = \\ = 
			\left\langle 
		f_{j_1p_1},  
		\left.\frac{\partial R_{j_1}(k,r)}{\partial r}\right|_{r=1}
			f_{j_1p_1}\right\rangle
		=
		 c_{j_1p_1}\left.\frac{\partial R_{j_1}(k,r)}{\partial r}\right|_{r=1}
		 \end{aligned}
	\end{equation}
	Combining \eqref{as_j}, (\ref{Lambda1}), (\ref{Lambda2}) and \eqref{5.28}, we obtain that
		\begin{equation}\label{l2sol2}
			\left\langle  f_{j_1p_1}, \left(\hat{\Phi}_1(E) - \hat{\Phi}_2(E)\right) f_{j_2p_2}\right\rangle \leq
			C(d) 2^{-j_1}||\psi_1 - \psi_2||_{\mathbb{L}^2(B(0,1))}.
		\end{equation}
	Combining (\ref{l2sol})	and (\ref{l2sol2}), we get (\ref{eq_3.2}) for $j_1\geq j_2$ and $E \neq 0$. 
	
  For $j_1 < j_2$ we use the fact that $\Phi_v^*(E) = \Phi_{\bar{v}}(\bar{E})$ in order to swap  $j_1$ and $j_2$, 
  where $\Phi_v^*$ denotes the adjoint operator to $\Phi_v$.
 Thus we complete the proof of Lemma \ref{Lemma_3.2} for the non-zero energy case.
  
  Estimate \eqref{eq_3.2} for the  zero energy case follows from Lemma 1 of \cite{Mandache2001}.
\section{Acknowledgements}
This work was fulfilled in the framework of research carried out under the
supervision of R.G. Novikov.

\noindent
{ {\bf M.I. Isaev}\\
Centre de Math\'ematiques Appliqu\'ees, Ecole Polytechnique,

91128 Palaiseau, France\\
Moscow Institute of Physics and Technology,

141700 Dolgoprudny, Russia\\
e-mail: \tt{isaev.m.i@gmail.com}}\\




\begin{thebibliography}{99}

\bibitem{Alessandrini1988}
 G. Alessandrini, 
 {\it Stable determination of conductivity by boundary measurements}, 
 Appl.Anal. 27, 1988, 153-172.

\bibitem{AV2005}
G. Alessandrini, S. Vassella,
{\it Lipschitz stability for the inverse conductivity problem},
Adv. in Appl. Math. 35, 2005, no.2, 207-241.





\bibitem{BK2012}
L. Beilina, M.V. Klibanov, 
{\it Approximate global convergence and adaptivity for coefficient inverse problems}, Springer (New York), 2012. 407 pp.

\bibitem{Buckhgeim2008}
 A. L. Buckhgeim, 
 {\it Recovering a potential from Cauchy data in the two-dimensional
case}, J. Inverse Ill-Posed Probl. 16, 2008, no. 1, 19-33.


\bibitem{Calderon1980} Calder\'on, A.P., {\it On an inverse boundary problem}, Seminar on Numerical Analysis
and its Applications to Continuum Physics, Soc. Brasiliera de Matematica, Rio de Janeiro, 1980, 61-73.

 

\bibitem{CR2003}
M. Di Cristo and L. Rondi
{\it Examples of exponential instability for inverse
inclusion and scattering problems}
{ Inverse Problems}. 19 (2003) 685Ц701.





\bibitem{Gelfand1954}
I.M. Gelfand, {\it Some problems of functional analysis and algebra}, Proceedings of the
International Congress of Mathematicians, Amsterdam, 1954, pp.253-276.









\bibitem{IsaevDtN}
M.I. Isaev,
{\it Exponential instability in the Gel'fand inverse 
problem on the energy intervals}, J. Inverse Ill-Posed Probl., Vol. 19(3), 2011, 453-473; 
e-print arXiv: 1012.2193. 

\bibitem{IN2012}
M.I. Isaev, R.G. Novikov
{\it Stability estimates for determination of potential
from the impedance boundary map}, e-print arXiv:1112.3728. 

\bibitem{IN2012+}
M.I. Isaev, R.G. Novikov
{\it Reconstruction of a potential
from the impedance boundary map}, e-print arXiv:1204.0076. 

\bibitem{IN2012++}
M.I. Isaev, R.G. Novikov
{\it Energy and regularity dependent stability estimates for the Gel'fand inverse problem in multidimensions}, 
e-print hal-00689636. 

\bibitem{Isakov2011}
V. Isakov, {\it Increasing stability for the  Schr\"odinger potential from the Dirichlet-to-Neumann map},
Discrete Contin. Dyn. Syst. Ser. S 4, 2011, no. 3, 631-640.


\bibitem{KV1985}
R. Kohn, M. Vogelius, {\it Determining conductivity by boundary measurements II}, Interior
results, Comm. Pure Appl. Math. 38, 1985, 643-667.




\bibitem{LR1986}
M.M. LavrentТev,  V.G. Romanov, S.P. Shishatskii, 
{\it Ill-posed problems of mathematical physics and analysis}, Translated from the Russian by J. R. Schulenberger.
Translation edited by Lev J. Leifman. Translations of Mathematical Monographs, 64.
American Mathematical Society, Providence, RI, 1986. vi+290 pp.

\bibitem{Liu1997}
L. Liu, 
{\it Stability Estimates for the Two-Dimensional Inverse Conductivity Problem},
Ph.D. thesis, Department of Mathematics, University of Rochester, New York, 1997.



\bibitem{Mandache2001}
 N. Mandache,
{ \it Exponential instability in an inverse problem for the Schr\"odinger equation}, 
 { Inverse Problems}. 17, 2001, 1435-1444.  



\bibitem{Nachman1996}
A. Nachman, {\it Global uniqueness for a two-dimensional inverse boundary value problem},
Ann. Math. 143, 1996, 71-96.

\bibitem{NUW2011}
S. Nagayasu, G. Uhlmann, J.-N. Wang, 
{\it Increasing stability in an inverse problem for the acoustic equation},
e-print arXiv:1110.5145



\bibitem{Novikov1988}
R.G. Novikov, {\it Multidimensional inverse spectral problem for the equation $-\Delta \psi + (v(x) - Eu(x))\psi = 0$}
Funkt. Anal. Prilozhen. 22(4), 1988, 11-22 (in Russian); Engl. Transl.  Funct. Anal. Appl. 22, 1988, 263-272.



\bibitem{Novikov1998}
R.G. Novikov, {\it Rapidly converging approximation in inverse quantum scattering in dimension
2}, Physics Letters A 238, 1998, 73-78.


 

\bibitem{Novikov2005+}
R.G. Novikov,
{ \it The $\bar\partial$-approach to approximate inverse scattering at fixed energy in three dimensions. }
 IMRP Int. Math. Res. Pap. 2005, no. 6, 287-349. 

\bibitem{Novikov2005}
R.G. Novikov, {\it Formulae and equations for finding scattering data from the Dirichlet-to-Neumann map with nonzero background potential}, 
Inverse Problems 21, 2005, 257-270.


\bibitem{Novikov2008}
R.G. Novikov, 
{ \it The $\bar\partial$-approach to monochromatic inverse scattering  in three dimensions,}
J. Geom. Anal 18, 2008, 612-631.





\bibitem{Novikov2011}
R.G. Novikov,
{\it New global stability estimates for the Gel'fand-Calderon inverse problem,} Inverse 
Problems 27, 2011, 015001(21pp); e-print arXiv:1002.0153.


\bibitem{NS2010} R. Novikov and M. Santacesaria, {\it A global stability estimate for 
the GelТfand- Calderon inverse problem in two dimensions}, J.Inverse 
Ill-Posed Probl., Volume 18, Issue 7, 2010, Pages 765-785; e-print arXiv: 1008.4888. 


\bibitem{NS2012}
R. Novikov and M. Santacesaria, 
{\it Monochromatic Reconstruction Algorithms for Two-dimensional Multi-channel Inverse Problems},
International Mathematics Research Notes, 2012, doi: 10.1093/imrn/rns025.



\bibitem{Rondi2006}
L. Rondi, {\it A remark on a paper by Alessandrini and Vessella}, Adv. in Appl. Math. 36 (1), 2006, 67-69.



 
\bibitem{S2011N}
M. Santacesaria, 
{\it New global stability estimates for the Calderon inverse problem in two dimensions},
 e-print: hal-00628403.
 
\bibitem{S2012} 
M. Santacesaria, {\it
Stability estimates for an inverse problem
for  the Schr\"odinger equation at negative
energy in two dimensions},
e-print: hal-00688457.
 
  
\bibitem{SU1987}
J. Sylvester and G. Uhlmann,  {\it A global uniqueness theorem for an inverse boundary
value problem}, Ann. of Math. 125, 1987, 153-169.




\bibitem{Weyl1912}
H. Weyl, {\it Das asymptotische Verteilungsgesetz der Eigenwerte linearer partieller Differentialgleichungen} (mit einer Anwendung auf die Theorie der Hohlraumstrahlung). Math. Ann. 71(4), 441-479, 1912.


\end{thebibliography}
\end{document}